\pdfoutput=1 

\documentclass[11pt]{amsart}
\usepackage{amsmath,amsthm,amssymb,amscd,eucal}        
\usepackage{graphicx}
\usepackage[all]{xy}
\usepackage[margin=1.3in]{geometry}

\begin{document}

\newcommand{\half}{{\textstyle{\frac{1}{2}}}}
\newcommand{\balpha}{{\boldsymbol{\alpha}}}
\newcommand{\B}{{\bf 1}}
\newcommand{\e}{{\bf e}}
\newcommand{\y}{{\bf r}}
\newcommand{\x}{{\bf x}} 
\newcommand{\tx}{{\bar{x}}}                   
\newcommand{\bn}{{\bf n}}

%
%

\newcommand{\aff}  {{\rm Aff}_+}                 
\newcommand{\Conf}{{\rm Config}}                 
\newcommand{\conf}{{\rm config}}                 
 
\newcommand{\PSL} {\Pj\Sl_2(\R)}                           
\newcommand{\PGL} {\Pj\Gl_2(\R)}                           
\newcommand{\PGLC} {\Pj\Gl_2(\C)}                          
\newcommand{\RP} {\R\Pj^1}                                 
\newcommand{\CP} {\C\Pj^1}                                 
\newcommand{\C} {{\mathbb C}}                              
\newcommand{\R} {{\mathbb R}}                              
\newcommand{\Z} {{\mathbb Z}}                              
\newcommand{\Pj} {{\mathbb P}}                             
\newcommand{\Sg} {{\mathbb S}}                             
\newcommand{\Gl} {{\rm Gl}}                                
\newcommand{\Sl} {{\rm Sl}}                                
\newcommand{\SO} {{\rm SO}}                                
\newcommand{\orb}{{\rm orb}}
\newcommand{\PGl}{{\rm PGl}} 
\newcommand{\SP}{{\rm SP}}
\newcommand{\Ot}{{\rm O}}
\newcommand{\D}{{\rm D}}
\newcommand{\reg}{{\rm reg}}

\newcommand{\cA}{{\mathcal A}}
\newcommand{\cB}{{\mathcal B}}
\newcommand{\cC}{{\mathcal C}}
\newcommand{\blambda}{{\boldsymbol{\lambda}}}
\newcommand{\bDelta}{{\boldsymbol{\Delta}}}
\newcommand{\bdelta}{{\boldsymbol{\delta}}}
\newcommand{\cQ}{{\mathcal Q}}
\newcommand{\dQ}{{\mathcal Q}^0}                               
\newcommand{\qQ}{{\widehat{\cQ}}}
\newcommand{\tQ}{\widetilde{\cQ}}
\newcommand{\sQ}{{\sf Q}}

\newcommand{\cell}{{\rm Cell}}
\newcommand{\cells}{{\rm Cells}}
\newcommand{\iso}{{\rm Iso}}

\newcommand{\metmap}{\mathfrak{d}}

\newcommand{\dis} {{\mathcal D}}                                 
\newcommand{\diagonal}{\triangle}                                

\newcommand{\oM} [1] {\ensuremath{{\mathcal M}_{0,#1}(\R)}}               
\newcommand{\M} [1] {\ensuremath{{\overline{\mathcal M}}{_{0, #1}(\R)}}}  
\newcommand{\OM} [1] {\ensuremath{{\overline{\mathcal M}}{_{0, #1}^{{\rm or}}(\R)}}}
\newcommand{\KM} [1] {\ensuremath{{\overline{\mathcal M}}{_{0, #1}^{{\rm kap}}(\R)}}}
\newcommand{\cM} [1] {\ensuremath{{\mathcal M}_{0, #1}}}                  
\newcommand{\CM} [1] {\ensuremath{{\overline{\mathcal M}}_{0, #1}}}       

\newcommand{\sM} [1] {\ensuremath{{\sf M}_{0,#1}(\R)}}
\newcommand{\tsM} [1]{\ensuremath{{\widetilde{\sf M}}{_{0, #1}(\R)}}}

\newcommand{\tsQ} [1] {{\widetilde{\sf Q}}}         

\newcommand{\iDelta}{{\overset{o}\Delta}}

\newcommand{\bT}{{\overline{\mathcal T}}}
\newcommand{\T}[1]{{\mathcal T}_#1}                                         
\newcommand{\pT}[1]{{\bf T}_{#1}(\R)}             
\newcommand{\BHV}[1]{{\rm BHV}_{#1}}                    
\newcommand{\iBHV}[1]{{\rm BHV}_{#1}^+}                 

\newcommand{\hide}[1]{}

\newcommand{\suchthat} {\:\: | \:\:}
\newcommand{\ore} {\ \ {\it or} \ \ }
\newcommand{\oand} {\ \ {\it and} \ \ }

%
%

\theoremstyle{plain}
\newtheorem{thm}{Theorem}
\newtheorem{prop}[thm]{Proposition}
\newtheorem{cor}[thm]{Corollary}
\newtheorem{lem}[thm]{Lemma}
\newtheorem{conj}[thm]{Conjecture}
\newtheorem*{thmhull}{Theorem 11}
\theoremstyle{definition}
\newtheorem*{defn}{Definition}
\newtheorem*{exmp}{Example}

\theoremstyle{remark}
\newtheorem*{rem}{Remark}
\newtheorem*{hnote}{Historical Note}
\newtheorem*{nota}{Notation}
\newtheorem*{ack}{Acknowledgments}
\numberwithin{equation}{section}


\title{Diagonalizing the Genome I: Navigation in Tree Spaces}

\author{Satyan L.\ Devadoss}
\address{S.\ Devadoss: Williams College, Williamstown, MA 01267}
\email{satyan.devadoss@williams.edu}

\author{Jack Morava}
\address{J.\ Morava: Johns Hopkins University, Baltimore, MD 21218}
\email{jack@math.jhu.edu}

\begin{abstract}
The orientable cover of the moduli space of real genus zero algebraic curves with marked points is a
compact aspherical manifold tiled by associahedra, which resolves the singularities of the space of 
phylogenetic trees. The resolution maps planar metric trees to their underlying abstract representatives,
collapsing and folding an explicit geometric decomposition of the moduli space into cubes. This decomposition
endows the resolving space with an interesting canonical pseudometric. \medskip

\noindent
The second part of this paper defines a related (stacky) resolution of a space of real quadratic forms, and suggests,
perhaps without much justification, that systems of oscillators parametrized by such objects may
provide useful models in genomics. 
\end{abstract}



\maketitle

\baselineskip=17pt

%
%
\section{Introduction} 
\label{s:intro}

A classical problem in computational biology is the construction of a phylogenetic tree from gene sequence alignments for $n$ 
species. Billera, Holmes, and Vogtmann \cite{bhv} construct a space $\BHV{n}$ of such metric trees, with nonpositive curvature; 
this makes geometric methods (geodesics, centroids) available. Unfortunately, their space is not a manifold but a cone over a 
relatively singular simplicial complex. This paper introduces a compact orientable hyperbolic manifold $\OM{n+1}$, tiled by $n!$ 
Stasheff associahedra $K_n$, which resolves [Theorem~\ref{t:moduli-bhv}] the singularities of $\BHV{n}$.

\begin{figure}[h]
\includegraphics[width=.4\textwidth]{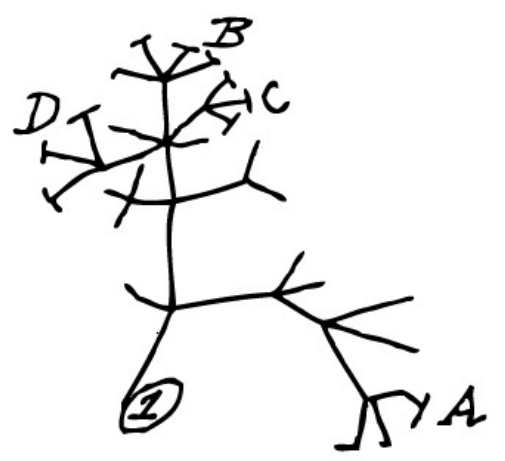}
\caption{An example from Darwin's notebook (1837).}
\label{f:darwin}
\end{figure}

Our manifold is the orientation cover of the moduli space \M{n+1} of real stable genus zero algebraic curves marked with 
distinct smooth points; these objects can also be interpreted as compactified moduli spaces of ordered 
configurations of points on the real projective line \cite{dev1}, but we understand them here as spaces of rooted metric 
trees with labeled leaves, as in Figure~\ref{f:darwin}. One of the objectives of this paper is to provide \M{n+1} with 
a pseudometric naturally associated to the geometry of the trees it parametrizes. 

Our spaces have some remarkable similarities \cite{mg} to the \emph{space forms} which captured the
imagination of mathematicians in the late nineteenth century: they are $K(\pi,1)$ manifolds with very pretty
tessellations. The work of Klein and other geometers soon found applications in physics and mechanics; we
hope our constructions may be useful for the study of big data in general, and biology in particular. 

\begin{ack}
We would like to thank Colin Adams, Mike Davis, Charles Epstein,
Sikimeti Mau, Amnon Ne'eman, and especially Jim Stasheff for continued interest and encouragement over many
years. SD also thanks Lior Pachter and the University of California, Berkeley, and Clemens Berger and the
Universit\'{e} de Nice for their hospitality during his sabbatical visits.  JM was partially supported by DARPA's 
FunBio program, and his part in this would have never been possible without Ben Mann's visionary work 
spearheading that project.  
\end{ack} 

\noindent
\emph{Organization.} Sections \ref{s:trees} and \ref{s:planar} of this paper summarize work of Billera, Holmes, and
Vogtmann on spaces of abstract metric trees, and the corresponding classical work of Boardman and Vogt on spaces of planar trees.  Sections \ref{s:config} through \ref{s:metric} use previous work of Devadoss to construct an \emph {orientable} moduli space of planar trees (a prerequisite for statistical applications), and Section \ref{s:maps} concerns the (equivariant) topology of the resolving map from the space of planar to abstract trees. 

The second part of this paper is a sketch, the expression of pious hope for future 
progress. It uses the constructions of Part I to define a topological stack of quadratic forms
with labeled eigenvalues, imagined as systems of interacting oscillators which maintain
their temporal integrity. This is related to the Coxeter complex of type $A$, the braid hyperplane arrangement.  The resulting geometric object has rather mild singularities; we hope it 
will be useful as a proving ground for evolutionary models.

%
%
\section{Spaces of abstract trees} 
\label{s:trees}

\begin{figure}[h]
\includegraphics{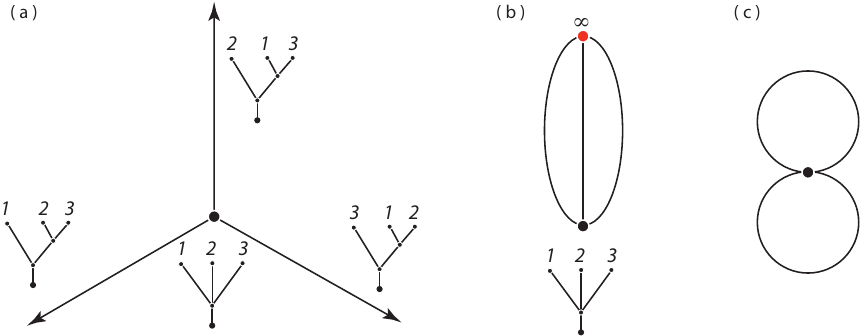}
\caption{The space (a) $\BHV{3}$, (b) $\iBHV{3}$, and (c) its homotopy equivalent.}
\label{f:bhv3}
\end{figure}

A \emph{metric tree} is a tree with a nonnegative weight assigned to each of its \emph{internal} edges
\cite{bd}. Billera, Holmes, and Vogtmann \cite{bhv} have constructed a space $\BHV{n}$ of isometry
classes of metric trees with $n$ labeled leaves (including a \emph{root}). When such a tree is binary, it has $n-2$
internal nodes, and hence $2n-2$ nodes in all. Diaconis and Holmes \cite{dh} have shown that such trees,
with all nodes labelled, are in bijection with the $(2n-3)!!$ possible (unordered) pairings of those nodes.
This enumeration defines coordinate patches for the space of such trees, parametrized by the weights of their
internal edges. The resulting $\BHV{n}$ space is contractible (in fact, a cone), formed by gluing $(2n-3)!!$ orthants
$\R_{\geq 0}^{n-2}$, one for each type of labeled binary tree.  As the weights go to zero, we get
degenerate trees on the boundaries of the orthants. Two boundary faces are identified when they contain the
same degenerate trees.  Figure~\ref{f:bhv3}(a) shows that $\BHV{3}$ consists of 3 rays $\R_{\geq 0}$ glued together at the origin.

Boardman's subspace $\T{n} \subset \BHV{n}$ of \emph{fully-grown} rooted trees \cite{boa} consists of trees with weights 
in $[0,1]$, having at least one edge of weight 1. Note that the subspace of $\BHV{n}$ of trees with weights that sum to 1 is 
also homeomorphic to $\T{n}$, by rescaling the maximal edge weight to 1. If a tree is not a star (or \emph{corolla}, if rooted)
with all internal weights zero, then the tree can be rescaled (by dividing its weights by their maximum) to be full-grown. 
This defines a homeomorphism with the (unreduced) smash product
\[
\BHV{n} \ \simeq \ (\R_{\geq 0},0) \wedge \T{n}  
\]
which sends a tree with edges of weight $\{w_k\}$ to the pair defined by $\max\{w_k\} \in \R_+$, together with a full-grown 
rescaled version of that tree; the star is sent to the cone point $0 \times \{\T{n}\}$ of $\BHV{n}$. In particular:

\begin{lem}
The one-point compactification $\iBHV{n}$ is homeomorphic to the (unreduced) suspension of $\T{n}$.
\end{lem}

For example, the subspace $\T{3}$ of $\BHV{3}$ is simply three points, one for each ray.  The one-point compactification 
$\iBHV{3}$ is the suspension of these three points, resulting in part (b) of the figure.  This is homotopy equivalent 
to the wedge of two circles, shown in (c). An early result by Vogtmann \cite{vog} related to outer automorphism groups 
generalizes this observation:

\begin{prop}
\label{p:wedge-sphere}
There is a homotopy equivalence
\begin{equation}
\psi \ : \ \iBHV{n} \ \longrightarrow \  \bigvee_{(n-1)!} S^{n-2} \;,
\label{e:rw}
\end{equation}
between $\iBHV{n}$ and a wedge of $(n-1)!$ spheres of dimension $n-2$. 
\end{prop}

\begin{figure}[h]
\includegraphics{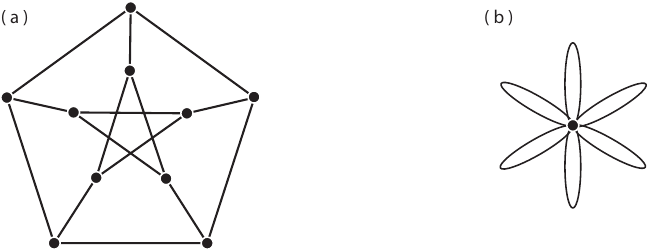}
\caption{The space (a) $\T{4}$ as the Peterson graph, where contracting certain edges results in (c) the wedges of six circles.}
\label{f:bhv4}
\end{figure}

\begin{exmp}
Consider the case of rooted trees with $n=4$ leaves.  The space $\BHV{4}$ consists of 15 quadrants $\R^2_{\geq 0}$ glued 
together.  The subspace $\T{4}$ is displayed in Figure~\ref{f:bhv4}(a) as the Peterson graph.  The 10 vertices correspond 
to rooted trees with four leaves with one internal edge of weight 0, and the 15 edges belong to the binary trees with no 
internal zero weights.  Part (b) shows the homotopy type of $\T{4}$, a wedge of six circles, obtained by contracting certain edges of (a).  Since $\iBHV{4}$ is the suspension of $\T{4}$, its homotopy type is wedge of six 2-spheres.
\end{exmp}

The spaces $\T{n}$ of trees have an important literature in representation theory and combinatorics \cite{rw}, also appearing 
in \cite{ss} as the \emph{tropical Grassmannian} ${\mathcal G}'''_{2,n+1}$. The symmetric group $\Sg_n$ acts on $\T{n}$ by 
permuting the labels of the leaves; this restricts to an $(\Sg_{n-1})$-equivariant homotopy equivalence 
\begin{equation}
\T{n} \; \sim \; S^{n-3} \wedge \Sg_{n-1}^+ \;.
\label{e:wedge}
\end{equation}
In fact, by dropping the distinction between root and leaf, the group action extends to $\Sg_{n+1}$. On the other hand, the 
suspension of $\iBHV{n}$ is $\Sg_n$-equivariantly homotopy-equivalent to a version of the geometric realization of the poset 
of partitions of $n$ \cite[\S 9]{chi}. The associated homology representations are fundamental in (among other things) 
the theory of operads. 

%
%
\section{Planar Trees and Cubes} 
\label{s:planar}
\subsection{}

An elegant polytope captures the structure of the space of planar rooted trees:

\begin{defn} 
The \emph{associahedron} is a convex  polytope of dimension $n-2$ whose face poset is isomorphic to that of bracketings of $n$ letters, ordered so $a \prec a'$ if $a$ is obtained from $a'$ by adding new brackets.
\end{defn}

\begin{figure}[h]
\includegraphics{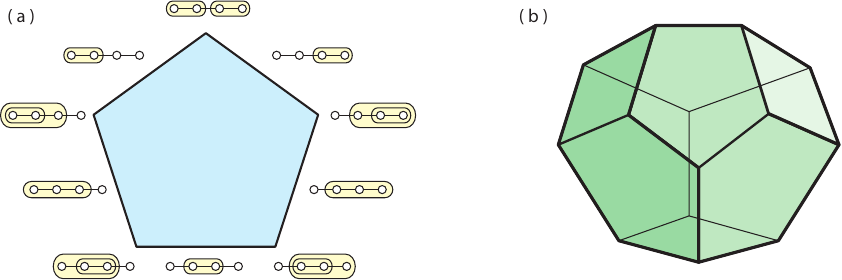}
\caption{Associahedra $K_4$ and $K_5$.}
\label{f:assoc}
\end{figure}

\begin{figure}[h]
\includegraphics{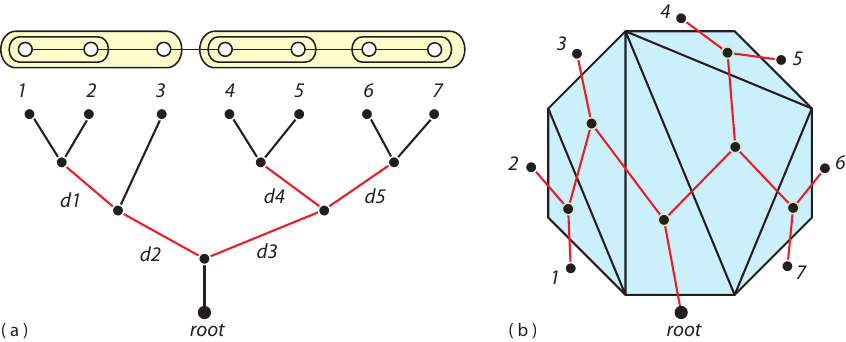}
\caption{(a) Bracketings and planar trees and (b) polygons with diagonals.}
\label{f:assoc-tree}
\end{figure}

\noindent
The  associahedron was constructed independently by Haiman (unpublished) and Lee \cite{lee}, though
Stasheff had defined the underlying abstract object twenty years previously, in his work on associativity in
homotopy theory \cite{sta}.  Figure~\ref{f:assoc}(a) shows the 2D associahedron $K_4$ with a labeling of
its faces, and (b) shows the 3D version $K_5$.
There are over a hundred combinatorial and geometric interpretations of the Catalan number
\begin{equation} \label{e:cubes}                            
C_{n-1} \ = \ \frac{1}{n}\binom{2n-2}{n-1} \ ,
\end{equation} 
which index the vertices of the associahedron $K_n$. Most important to us is the relationship between bracketings of $n$ letters,
which capture collisions of points on the line, and rooted planar binary trees. Figure~\ref{f:assoc-tree}(a)
illustrates the bijection with rooted trees, while (b) shows the relationship to polygons with diagonals.

\begin{thm}
\label{t:cubes}
The \emph{W-construction} of Boardman and Vogt defines a decomposition of $K_n$ into 
$C_{n-1}$ $(n-2)$-dimensional cubes, which assigns it a canonical pseudometric.
\end{thm}

\begin{proof}
A tree embedded in the plane inherits an orientation. To define coordinate patches on the space of 
(isometry classes of) \emph{planar} rooted metric trees with $n$ leaves, labeled $\{1,2, \ldots, n\}$ 
clockwise, we assign weights to the internal edges defined by a plane embedding.  In particular, terminal edges are 
to be understood as having infinite weight. Such data defines a parametrization of the space of such planar trees by a 
Catalan number $C_{n-1}$ of orthants $\R^{n-2}_{\geq 0}$, which glue together as in \cite[\S 1.19]{bv}.  Completing 
the orthants to cubes $[0,\infty]^{n-2}$ \cite{bhv} defines the manifold-with-corners $K_n$ of planar labeled trees, whose faces 
parametrize degenerate trees, with at least one internal edge length equal to $\infty$.
\end{proof}

\noindent We will describe the resulting metric in more detail in Section~\ref{s:metric}. 

\begin{figure}[h]
\includegraphics{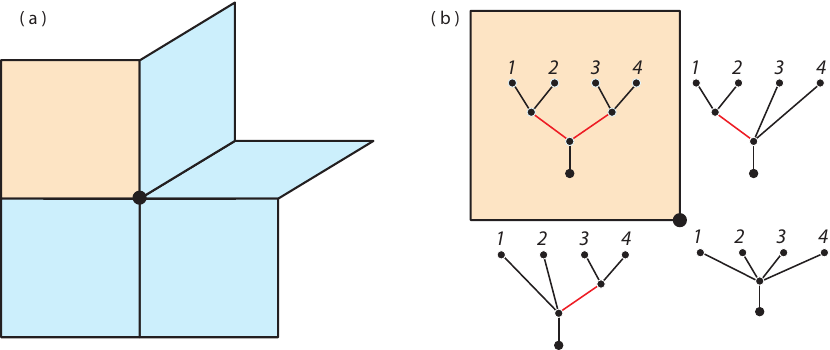}
\caption{The piecewise-Euclidean geometry of $K_4$.}
\label{f:bv-trees}
\end{figure}

\begin{exmp}
Figure~\ref{f:bv-trees}(a) shows the construction for $n=4$ by gluing five quadrants (drawn as squares), along with a 
labeling of one of the squares in (b). In general, all the orthants will share a common vertex defined by the star.
\end{exmp}

\subsection{}

\begin{figure}[h]
\includegraphics{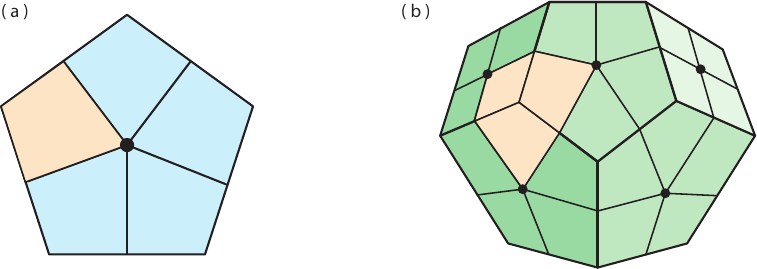}
\caption{The cubical decomposition of associahedra.}
\label{f:bv-cubes}                                                         
\end{figure}

The Boardman-Vogt construction is a decomposition of the associahedron into cubes, exemplifying the
general cubic barycentric subdivision of a polytope. Figure~\ref{f:bv-cubes} shows the decomposition 
for the examples from Figure~\ref{f:assoc}, where one cube from each complex has been highlighted. Each
vertex of $K_n$, and hence each cube of the complex, is associated to a planar binary tree with $n$ leaves;
there is thus a Catalan number of cubes in the decomposition.

\begin{figure}[h]
\includegraphics{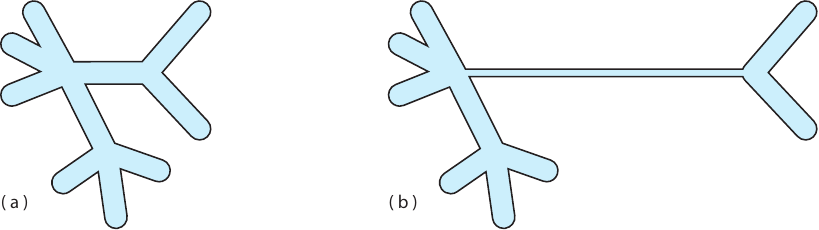}
\caption{A ribbon tree as the weight of an internal edge is altered.}
\label{f:ribbon}
\end{figure}

These planar trees can be viewed as \emph{ribbon trees}, with a natural cyclic order around each vertex
given by the planar embedding; see Figure~\ref{f:ribbon}(a).  As the weight of an internal edge increases,
we think of the thickness of the ribbon as decreasing, as in part (b).  As the length approaches $\infty$, the
edge becomes a ribbonless segment. 
Faces of $K_n$ are canonically products $\prod K_{\nu_i-1}$ of lower-dimensional associahedra, indexed by
valency $\nu_i$ of the internal nodes of rooted trees with $(n+1)$ leaves. 

The associahedra have some important symmetries. Relabeling the leaves of a tree in clockwise order (that is,
turning over the piece of paper on which the tree is drawn) defines an involution $\varepsilon: K_n \to
K_n$. A related construction interprets the root as a \emph{zeroth} leaf, shifting the resulting indexing by a
one-step clockwise rotation $r$.  This defines an action of cyclic group $\Z_{n+1}$ on $K_n$ which,
together with the involutions above, extends to an action of the dihedral group $D_{n+1}$, equivalent to that defined by Lee \cite{lee}. This points to the usefulness of trees with arbitrarily labeled leaves, indexed by elements of $\Sg_n$ (if the root is distinguished) or by $\Sg_{n+1}$ (if not).  We now turn to considering a space of such labeled planar trees, tiled by the associahedron polytope.

%
%
\section{Configuration spaces}
\label{s:config}
\subsection{}

The moduli problem for algebraic curves has been a central problem in mathematics since Riemann.  In the 1970s 
it was solved over the integers $\Z$ by Deligne, Mumford, Knudsen \cite{kn} and others. A very
special case of their results constructs a moduli space for \emph{real} algebraic curves of genus zero
marked with $n \geq 3$ distinct smooth points.  That solution can be regarded as a good compactification \M{n+1}
of the space
$$\oM{n+1} \ = \ \Conf^{n+1}(\RP)/\PGL$$
of $n+1$ distinct particles on the real projective line. If we think of $\PGL$ as fixing three of the $n+1$ points 
$\{0, 1, 2, \ldots, n-1, \infty\}$ on $\RP$, the space \oM{n+1} can be reinterpreted as the $(n-3)$-torus $(\RP)^{n-2}$ with 
the subspace 
\begin{equation}
\Delta = \{(x_1, \ldots , x_{n-2}) \in (\RP)^{n-2} \suchthat x_i = x_j \ \ {\rm or} \ \ x_i = 0, 1, \infty\}
\label{e:delta}
\end{equation}
removed.   There is a stratification of $\Delta$ having
\begin{equation}
3 \binom{n-2}{k}   \ + \  \binom{n-2}{k-1}
\label{e:diag}
\end{equation}
cells of dimension $k$, corresponding to where $n-k-1$ particles collide.  The first term of \eqref{e:diag} counts collisions 
of free particles $x_i$ with one fixed particle $\{0,1,\infty\}$, whereas the second term counts collisions of only free particles.

Constructing the compactification of the real moduli space $\oM{n+1}$ uses the algebro-geometric notion of a \emph{blowup}.
For subspace $X$ of manifold $Y$, the \emph{blowup} of $Y$ along $X$ is obtained by first removing $X$ and replacing it with 
the sphere bundle associated to the normal bundle of $X \subset Y$. We then \emph{projectify} the bundle.

\begin{thm} \cite[\S 4]{dev1}
\label{t:tiling}
The Deligne-Mumford-Knudsen compactification of \oM{n+1} is obtained by performing iterated blow-ups along the cells of 
$\Delta$ in increasing order of dimension. The resulting space \M{n+1} is tiled by $n!/2$ copies of the associahedron $K_n$.
\end{thm}

It is easy to see that \M{3} is a point.  The manifold \M{4} is homeomorphic to a circle, with the {\em cross-ratio} serving 
as the homeomorphism; it is tiled by three $K_3$ line segments, glued together to form a triangle.

\begin{figure}[h]
\includegraphics{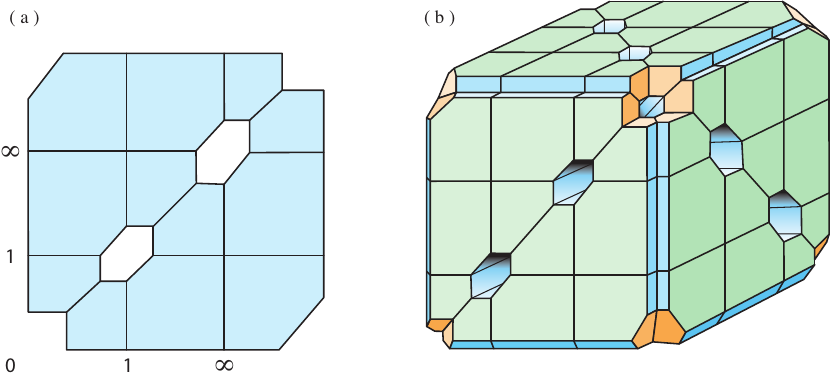}
\caption{(a) \M{5} and (b) \M{6} as blowups of tori.}
\label{f:m0506}
\end{figure}

\begin{exmp} 
An illustration of \M{5} from the construction in Theorem~\ref{t:tiling} appears in Figure~\ref{f:m0506}(a).  The diagonal 
$\Delta$ becomes the set of three points $\{x_1=x_2 = 0,1,\infty\}$; blowing up along these points yields \M{5} as the connected 
sum of a torus with three real projective planes, tiled by 12 associahedra $K_4$ from Figure~\ref{f:assoc}(a).  Indeed, for 
all $n \geq 4$, the moduli spaces \M{n+1} are nonorientable due to the real blowups along $\Delta$.
\end{exmp}

\begin{exmp}
Figure~\ref{f:m0506}(b) shows \M{6} as the blow-up of a 3-torus. The set $\Delta$ associated to \M{6} has three points 
\{$x_1=x_2=x_3=0,1,\infty$\} and  ten lines \{$x_i=x_j=0,1,\infty$\} and \{$x_1=x_2=x_3$\}.  The points are first blown up, 
each producing a connected sum with $\R\Pj^3$.  The lines are then blown up, corresponding to the hexagonal prisms, nine 
cutting through the faces, and the tenth (hidden) running through the torus from the bottom left to the top right corner.   
The moduli space \M{6} is tiled by 60 associahedra $K_5$ from Figure~\ref{f:assoc}(b). 
\end{exmp}

\subsection{}

Theorem~\ref{t:tiling} can be fleshed out combinatorially, by describing \M{n+1} as glued together from $n!/2$ associahedra;  
further details can be found in \cite[\S 3]{dev1}. Consider an $(n+1)$-gon\footnote{It can be helpful to think of 
these polygons as hyperbolic, with their vertices on the boundary of the Poincar\'e disk, and geodesic edges, as in Figure~\ref{f:hyperbolic} below.}  
with its sides labeled using $n$ labels $\{0, 1, 2, \ldots, n-1, \infty\}$ such that the the labels $\{0,1, \infty\}$ appear 
in clockwise (CW) order.  Up to rotation, this results in $n!/2$ possible labelings, each corresponding to an associahedral 
chamber of \M{n+1}.  Call these \emph{CW-labeled} polygons.  The faces of the associahedra are obtained by considering 
(nonintersecting) diagonals of these CW-labeled polygons, which, as in Figure~\ref{f:assoc-tree}, is analogous to resolving a tree.

\begin{defn}
A diagonal $d$ of a labeled $(n+1)$-gon $P$ partitions it into two smaller subpolygons.  The subpolygon labeled with at most 
one element of $\{0,1,\infty \}$ is called the \emph{lightgon}, and the complementary subpolygon the \emph{heavygon}.
\end{defn}

\begin{figure} [h]
\includegraphics[width=\textwidth]{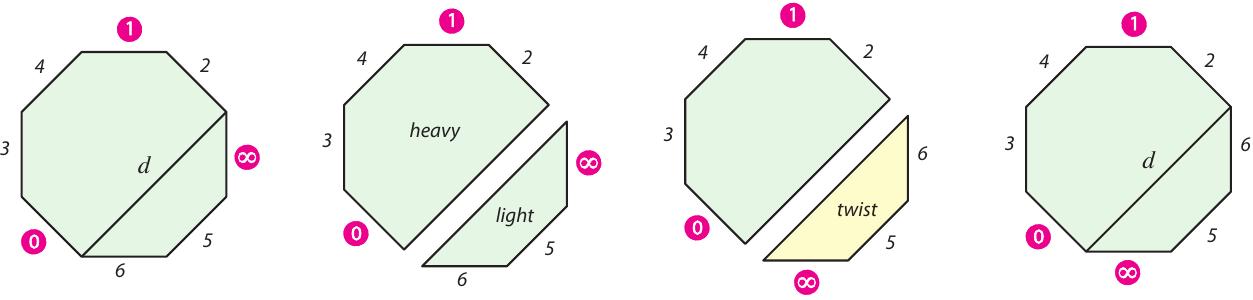}
\caption{Twisting along diagonal $d$.}
\label{f:dtwist}
\end{figure}

\begin{defn}
A \emph{twist} along diagonal $d$ of $P$ is the labeled polygon obtained by `breaking' $P$ along $d$ into two subpolygons, 
`reflecting' its lightgon, and `gluing' them back (Figure~\ref{f:dtwist}).  Notice that if a polygon is CW-labeled, it remains 
CW-labeled after a twist.
\end{defn}

\begin{thm} \cite[\S 3]{dev1}
\label{t:orig}
Two faces of $K_n$ associahedra, represented by CW-labeled polygons $P_1$ and $P_2$, are identified in \M{n+1} if twisting 
along certain diagonals of $P_1$ yields $P_2$.
\end{thm}

\begin{figure}[h]
\includegraphics{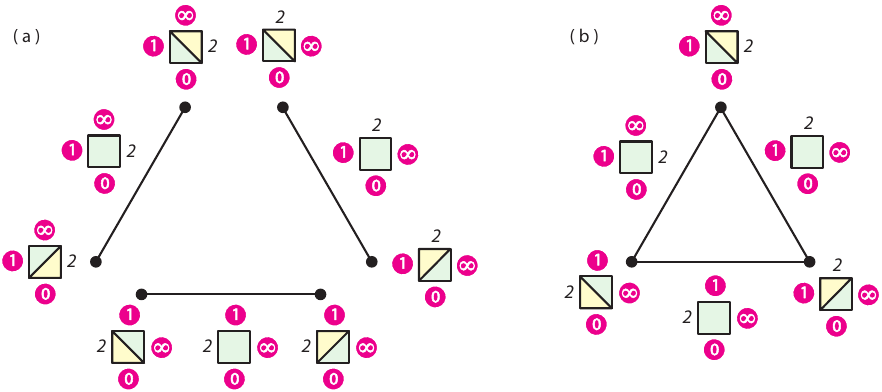}
\caption{Gluing three associahedra $K_3$ resulting in \M{4}.}
\label{f:gluing}
\end{figure}

\begin{exmp}
Figure~\ref{f:gluing} shows an example for \M{4}.  Part (a) displays the associahedron $K_3$ (an interval) with three possible 
labelings of the $4$-gon.  The boundaries of each associahedron correspond to drawing diagonals of the polygon, with the 
lightgons shaded in yellow.  Part (b) shows the gluing of the three associahedra under twisting, resulting in the triangle \M{4}.
\end{exmp}

The topology behind the combinatorics follows from the blowups, with details given in \cite[\S 5]{dev1}:  The $n+1$ particles 
colliding on $\RP$ correspond to the labels of the $(n+1)$-gon, and keeping track of such collisions along $\Delta$ is captured 
by the diagonals.  Notably, the labels along the lightgon-side of a polygon display the particles which have collided.  The 
twisting operation corresponds to the antipodal identification across the normal bundle due to the blow-up.

%
%
\section{The Orientable Cover}
\label{s:orient}
\subsection{}

We are now in position to describe and construct the orientable double cover \OM{n+1}.  Consider $n!$ copies of associahedra 
$K_n$, doubling the number of polytopes, represented by \emph{all} possible labelings of $(n+1)$-gons using $\{0, 1, 2, \ldots, 
n-1, \infty\}$.  Notice that half of them will be CW-labeled and the other half counterclockwise CCW-labeled.  If we only perform 
twists along diagonals of the polygons, then we would produce two copies of \M{n+1}, one with CW-labeling and another with 
CCW-labeling of the chambers.   Thus, another operation is needed to form the connected orientable cover.

\begin{figure} [h]
\includegraphics[width=\textwidth]{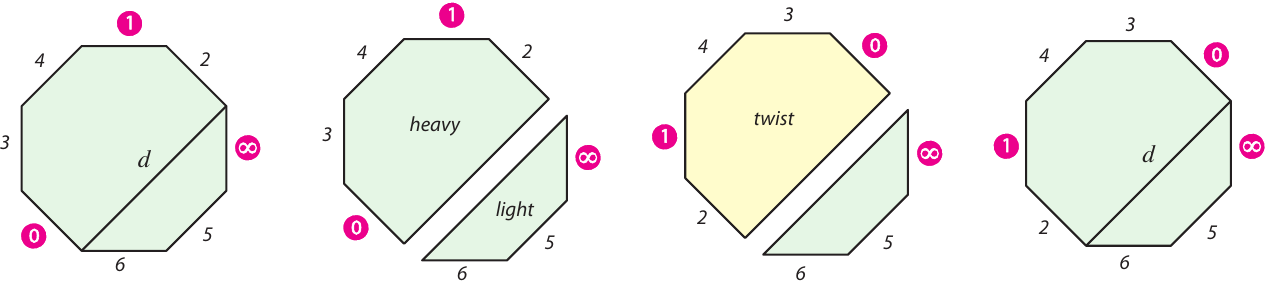}
\caption{Cotwisting along diagonal $d$.}
\label{f:cotwist}
\end{figure}

\begin{defn}
An \emph{cotwist} along diagonal $d$ of $P$ is the labeled polygon obtained first by `breaking' $P$ along $d$ into two 
subpolygons.  If the lightgon has an even number of \emph{labels},
reflect the lightgon and glue.  Otherwise, reflect the heavygon and glue (Figure~\ref{f:cotwist}).   
\end{defn}

\begin{thm}
\label{t:odc}
The orientable  cover \OM{n+1} is tiled by $n!$ copies of the associahedron $K_n$, identified with 
all possible labelings of an $(n+1)$-gon using $\{0, 1, 2, \ldots, n-1, \infty\}$.  
Two faces of associahedra, represented by labeled polygons $P_1$ and $P_2$, are identified if cotwisting along certain diagonals of $P_1$ yields $P_2$. 

Moreover, the symmetric group $\Sg_{n+1}$ acts smoothly on this cover by permuting labels, where the degree of $\sigma \in \Sg_{n+1}$  
is its sign as a permutation.
\end{thm}

\begin{proof}
Since the moduli space \M{n+1} is constructed from an orientable torus $(\RP)^{n-2}$, its nonorientability comes only from the
blowups. The blowup along a $k$-dimensional cell $X \subset \Delta$ of Eq.~\eqref{e:delta} sews in the real projective space 
$\R\Pj^{n-k-2}$ along $X$.  Such a cell corresponds to a collision of $n-k-1$ particles, represented by a polygon with a diagonal 
having $n-k-1$ labels on its lightgon.  Since $\R\Pj^i$ is orientable only for odd $i$, the blowup along a cell does not change 
orientability precisely when it has an even number of labels on its lightgons.

To construct the orientable cover, consider two copies of the torus $(\RP)^{n-2}$, as the configuration space of $n$ particles 
$\{0, 1, 2, \ldots, n-1, \infty\}$ on $\RP$ with labels $\{0,1, \infty\}$ fixed.  One copy has $\{0,1, \infty\}$ appearing in 
CW order, whereas the other in CCW order.  Recall that blowups along a cell $X$ are done in two stages, first replacing $X$ with 
the sphere bundle associated to the normal bundle; then projectify the sphere bundle.  For both tori, perform partial 
blowups (ie, just replace with the sphere bundle itself)  with along the cells of $\Delta$ in increasing order of dimension.  
By Theorem~\ref{t:tiling} this now produces two (disconnected) spaces composed of $n!$ associahedra $K_n$.

We glue these associahedra to form the connected, orientable double cover as follows:  First notice that reflecting the heavygon 
switches the parity of a polygon, between CW and CCW labeling, whereas reflecting the lightgon preserves the parity. 
If a facet of $K_n$ is a lightgon with even labels, reflect the lightgon.  This corresponds to gluing the facet antipodally across 
the sphere bundle of the same torus (a classic blowup).  By the discussion earlier, this blowup preserves orientability.  

On the other hand, if the lightgon has odd labels, twist along the heavygon.  This corresponds to gluing the facet to its 
antipodal counterpart on the other torus.   This gluing doubles the M\"obius strip, producing an orientable manifold.
The operation described above is is exactly embodied by the cotwist.
\end{proof}

\begin{rem}
Since \M{3} and \M{4} are orientable manifolds, the construction outlined in Theorem~\ref{t:odc} for these two cases results
in the trivial disconnected double cover.
\end{rem}

Kapranov \cite{kap} created an elegant \emph{non}orientable double cover of \M{n+1}, based on blowups of the braid hyperplane 
arrangement.  The $n!$ associahedra tiling this space are glued based on \emph{marked} twists \cite[\S 4]{dev1}: 
break a polygon into two pieces, reflect the side not containing the $\infty$  label, and glue.

\subsection{}

Quite a lot is now known about the topology of \M{n+1}, which extends to the orientable cover.  From a homotopy viewpoint, their 
tessellations imply their negative curvature in the sense of Alexandrov, and hence that they are $K(\pi,1)$ manifolds, whose 
fundamental groups are the twisted right-angled Coxeter groups of the associahedra \cite{djs}. Generators and relations 
\cite[\S 3.1]{hk} are known for these groups,  which have some striking analogies with that of the pure braid groups. 

The integral homology $H_*(\M{n+1}, \Z)$ is calculated in \cite{rains}, isomorphic to the direct sum of the cohomology of the 
relevant ideals in the poset of odd set partitions \cite[\S 4]{val}.  In particular, the  homology groups have only 2-torsion.  
The Poincar\'{e} polynomial of \M{n+1} is calculated in \cite[\S 2]{ehkr} as
\[
P_{n+1}(t) \ =  \prod_{0 \leq k < (n-2)/2} (1+ (n-2k-2)^2t) \, .
\]
This specializes to the Euler characteristic of \M{n+1}, calculated in \cite[\S 3]{dev1}, which yields
\begin{equation}
\chi (\OM{n+1}) \ = \ 
\begin{cases}
0 & n \text{ odd}\\
(-2)^{\frac{n-2}{2}}(n-1)!!(n-3)!! & n \text{ even.}
\end{cases}
\label{e:euler}
\end{equation}
for the double cover. 

\begin{figure} [h]
\includegraphics{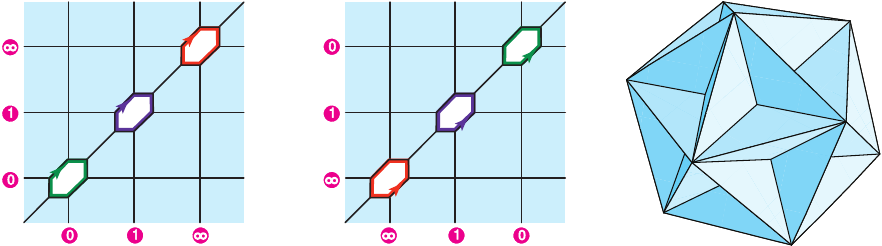}
\caption{(a) The CW and (b) CCW labeled spaces \M{5}, along with (c) the great dodecahedron appearing as the genus four \OM{5}.}
\label{f:gluing05}
\end{figure}

\begin{exmp}
The orientation cover \OM{5} can be constructed by considering two copies of \M{5} with CW and CCW labeling of their polygons, as displayed in Figure~\ref{f:gluing05}(a) and (b).  The blowups of the three 0-cells of $\Delta$ result in (nonorientable) crosscaps.  The orientable double cover is obtained by gluing the ``holes'' of the CW labeling to those of the CCW labeling, as given by the colored arrows, based on the cotwist operation.
The resulting manifold is orientable, tiled by $24$ associahedra $K_4$ pentagons.  Eq.~\eqref{e:euler} shows the Euler characteristic to be $-6$, forming a surface of genus four.  An embedding of this space with high symmetry is provided in \cite{as}, based on a pentagonal decomposition of Kepler's \emph{great dodecahedron} shown in Figure~\ref{f:gluing05}(c).
\end{exmp}

\begin{figure} [h]
\includegraphics{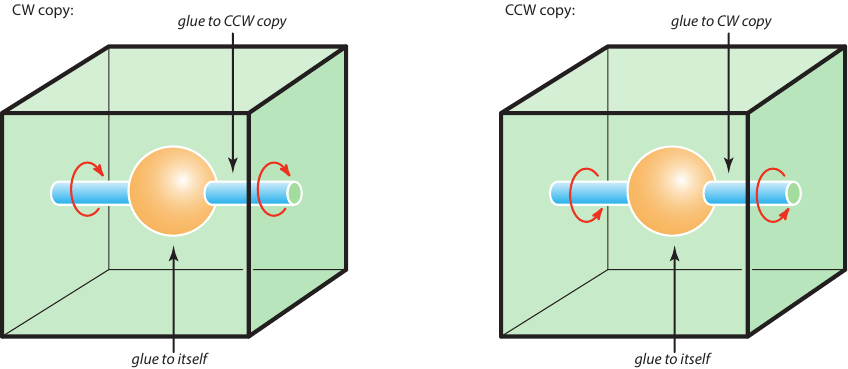}
\caption{Gluing by cotwisting to form \OM{6}.}
\label{f:gluing06}
\end{figure}

\begin{exmp}
Figure~\ref{f:m0506}(b) shows \M{6} obtained from the 3-torus, after blow-ups along three points and ten lines of $\Delta$.  Figure~\ref{f:gluing06} gives a closer view into a part of \OM{6}.  Since the blowup of each point results in sewing in an orientable $\R\Pj^3$ (in orange), each facet along the blowups glues antipodally to the same tori.  However, blowups along the ten lines result in gluing of $\R\Pj^2$, which introduce nonorientability along the tubes (in blue).  Thus, instead of an antipodal gluing across the normal bundle, we glue across to its antipodal counterpart on the other torus.  The resulting \OM{6} is an orientable 3-manifold tiled by 120  associahedra.
\end{exmp}

%
%
\section{Metrics on Moduli Spaces}
\label{s:metric}
\subsection{}

We can be more explicit about the piecewise Euclidean pseudometric on \OM{n+1} defined in Theorem~\ref{t:cubes}.
Let $\conf^n(\R)$ be the configuration space of $n$ distinct particles on $\R$ with the \emph{fixed linear ordering} from $v_1$ 
to $v_n$. Let 
\[
\aff \; := \: \{ x \mapsto ax + b \suchthat a > 0, \ b \in \R \}
\]
be the group of affine transformations of $\R$ generated by translations and positive dilations.   The action of $\aff$ 
on $\conf^n(\R)$ translates the left most particle $v_1$ to $0 \in \R$ and scales the remaining particles so that the maximum 
distance between any adjacent pairs $v_k$ and $v_{k+1}$ is set to $1$.

For a any rooted metric $T$, there is a unique path $\rho_k$ from the $k$-th leaf to its root. Let the weight $W_k$ be 
the sum of the lengths of the internal edges shared by the paths $\rho_k$ and $\rho_{k+1}$. Define
\[
v_k \ = \ \sum_{0 < i < k} \, e^{-W_i} \hspace{.5in} \textrm{for $k \in \{1, \dots, n-1\}$} \;,
\]
and let $v_1 = 0$.  The following relates configuration spaces to spaces of planar trees.

\begin{lem}
\label{l:metric}
The continuous embedding
\begin{equation} \label{e:metric}
\metmap \ : T \mapsto \{ 0 = v_1 < v_2 < \dots < v_n \ (< v_{n+1} = \infty) \} \ : \overset{o}{K_n} \ \longrightarrow \ 
\conf^n(\R)/\aff 
\end{equation}
identifies the interior of the associahedron with the space of configurations of distinct points on the line. 
\end{lem}

\begin{proof}
The map $\metmap$ takes vectors $(d_1, ..., d_{n-2})$ corresponding to interior edge lengths of a rooted tree with $n$ leaves 
to scaled ordered configurations $(v_1=0 < v_2 < .. <v_n)$ on the real line.  
For any configuration in the image, using the action of $\aff$, we can assume that there is always one distance between two 
points (say) $v_k$ and $v_{k+1}$ which is exactly one; this is the place where there is no overlap between $\rho_k$ and 
$\rho_{k+1}$, where they share no edges to the root.  Thus $e^0=1$, 
allowing us to partition the images into two separate 
groups, splitting along $v_{k-1}$ and $v_k$.  Induction for each group again, based on the lengths of the paths 
to the root, proves injection.

To reconstruct the tree from a configuration, start by partitioning at the places (maybe more than one) where $v_k$ and $v_{k+1}$ 
differ by exactly one.  This forms a first-level tree, with a root connected to each set of partitions.  Iterating this, where 
each partition set is rescaled to distance 1 and then repartitioned, forms subtrees. Stacking these trees (using operadic 
grafting) produces the reconstructed tree.
\end{proof}

\begin{exmp}
Figure~\ref{f:assoc-tree}(a) shows an example with seven leaves having five (colored) internal edges.  Here we have $W_1 = d_1 + d_2$, $W_2 = d_2$, $W_3 = 0$, $W_4 = d_3 + d_4$, $W_5 = d_3$ and $W_6 = d_3 + d_5$.   Based on this, Figure~\ref{f:metric-tree} gives the distances between the seven particles $v_i$ on $\R$, where $v_1$ is fixed at zero.
\begin{figure}[h]
\includegraphics{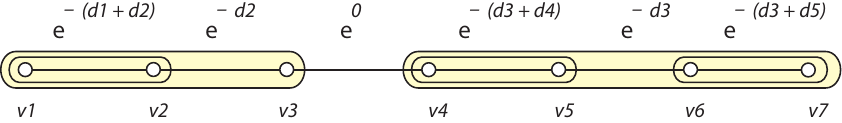}
\caption{The map $\metmap$ showing distances between particles on $\R$.}
\label{f:metric-tree}
\end{figure}
\end{exmp}

\begin{rem}
When the tree is a star (no internal edges) with $n$ leaves, then $W_k = 0$ for all $k$.
Under the map $\metmap$, the star tree is mapped to the configuration $(0, 1, 2, \cdots, n-1)$ of points in $\R$. 
On the other hand, if $W_k$ tends to $\infty$, then $(v_{k+1} - v_k)$ approaches $0$, and the corresponding points in the
configuration space collide. 
\end{rem}

\subsection{}

Lemma~\ref{l:metric} leads to one of our main results:

\begin{thm}
The tiling of \OM{n+1} by associahedra defines a psedometric extending that of Theorem 3.
\end{thm}

\begin{proof}
From \cite[Theorem 14]{dev2}, the Fulton-MacPherson compactification of $\conf^n(\R)/\aff$ yields an associahedral chamber 
of \OM{n+1}. The map $\metmap$ naturally extends to 
an embedding of the boundary of $K_n$ into this compactification, where the boundary of $K_n$ allows weights of certain tree 
edges to be infinite.  
\end{proof}

\begin{exmp}
Consider the one-cell of $\M{4}$ defined by configurations $0 < \rho < 1 <
\infty$ on the real projective line: the $W$-construction decomposes it into two one-cubes, the
first corresponding to (hyperbolic) trees as in Figure~\ref{f:hyperbolic}(a), the second as in part (b). In
either case there is one internal edge, of hyperbolic length $\lambda \in (0,\infty)$.

\begin{figure}[h]
\includegraphics{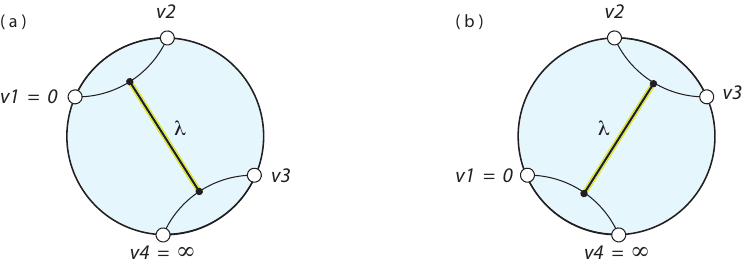}
\caption{Distances based on hyperbolic trees.}
\label{f:hyperbolic}
\end{figure}

In the first case, $\lambda \to \infty$ as $v_2 \to v_1$, and we have $W_1 = \lambda$ and  $W_2 =
0$.  Hence $v_1 = 0$,  $v_2 = e^{-\lambda}$,  and $v_3 = 1 + e^{-\lambda}$. In the second, $\lambda
\to \infty$ as $v_2 \to v_3$, and we have $W_1 = 0$ and $W_2 = \lambda$.  Hence $v_1 = 0$, $v_2 =
1$, and $v_3 = 1 + e^{-\lambda}$.
The normalization
\[
[z_0:z_1:z_2:z_3] \; := \; \frac{z_0 - z_1}{z_0 - z_2} \cdot \frac{z_3 - z_2}{z_3 - z_1}
\]
of the cross-ratio specializes to $(z_0 - z_1)(z_0 - z_2)^{-1}$ as $z_3 \to \infty$, and further to
$z_1z_2^{-1}$ as $z_0 \to 0$; thus
\[
[0:\rho:1:\infty] \; = \; \rho \;.
\]
If we think of $K_3 = [-\infty,0] \; \cup \; [0,+\infty]$ as partitioned into cubes, the first corresponding to 
trees with an internal leg of length $-\lambda$ (and the second to length $+\lambda$), then our metric
\[
\metmap(\lambda) = (1 + e^{-\lambda})^{-1} \ : \ K_3 = [-\infty,+\infty] \ \to \ [0,1] \subset \M{4}
\]
is the logistic function. 
\end{exmp}

\subsection{}

Algebraic geometers construct $\overline{\mathcal{M}}_{0,n}$ as the scheme which solves a certain moduli
problem. They then show it admits a projective embedding, and thus a nice metric; but there are
many such embeddings, corresponding to elements of a certain cone in its second cohomology, and there 
seems to be no obvious choice for a canonical element. The $\Sg_n$-equivariant embedding
\[
\M{n} \; \to \; \prod_{n \choose 4} \ \M{4}
\]
(the product taken over choices of four points out of $n$) seems to be an interesting candidate
for a canonical embedding, and hence for a canonical metric.

Similarly, the logarithmic derivative
\[
\alpha(z_0:z_1:z_2:z_3) \ := \ d\log[z_0:z_1:z_2:z_3] \ = \ A_{01} - A_{02} + A_{23} - A_{13}
\]
of the cross-ratio, where
\[
A_{ik} \ := \ (z_i - z_k)^{-1}(dz_i - dz_k) \ \in \ \Omega^1(\Conf^n(\C))
\]
is Arnol'd's one-form, defines an interesting volume form
\[
\bigwedge_{1 \leq i \leq n-2} \; \alpha(0:1:z_i:\infty) \ \in \ \Omega^{n-2} \; \OM{n+1} \ \,
\]
pulled back from $\otimes^{n-2} \; \Omega^1 \; \M{4}$. For example, the integral of 
\[
\alpha(1:0:\infty:s) \wedge \alpha(t:1:\infty:0)
\]
over $\overset{o}{K}{_4} = \{ 0 < s < t < 1 < \infty \}$ is the period
\[
\int_{0 < t < s < 1} s^{-1}ds \ \wedge \ (1-t)^{-1}dt \ = \ \zeta(2) \ = \ -\frac{(2\pi i)^2}{24}
\]
of a certain motive \cite[\S 25]{andre} studied by Goncharov and Manin \cite{gm}. 

%
%
\section{Maps between Tree Spaces}  \label{s:maps}
\subsection{}

We now provide the resolution maps from planar metric trees to their underlying abstract representatives,
collapsing and folding an explicit geometric decomposition of the moduli space into cubes. 

\begin{thm}
\label{t:moduli-bhv}
There is a smooth $\Sg_n$-equivariant blowup 
$$\Theta \ : \  \OM{n+1} \ \to \ \iBHV{n} \; $$
of the compactified Billera-Holmes-Vogtmann space by a compact smooth orientable aspherical manifold.  This map is 
generically $2^{n-1} : 1$, represented by an origami folding of cubes.  
\end{thm}

\begin{proof}
Partition each associahedron of $\OM{n+1}$ into $C_{n-1}$ number of cubes, according to Theorem~\ref{t:cubes}, one for each 
vertex.  Figure~\ref{f:map-bhv}(a) shows the cubical decomposition of four associahedra $K_4$ which will glue together, providing 
a local picture of \OM{5}.   Since \OM{n+1} is a right-angled \cite{djs} manifold, $2^{n-2}$ associahedra meet at each vertex;  
part (b) displays four small cubes meeting at a vertex of \OM{5}, forming a larger cube.  

\begin{figure}[h]
\includegraphics{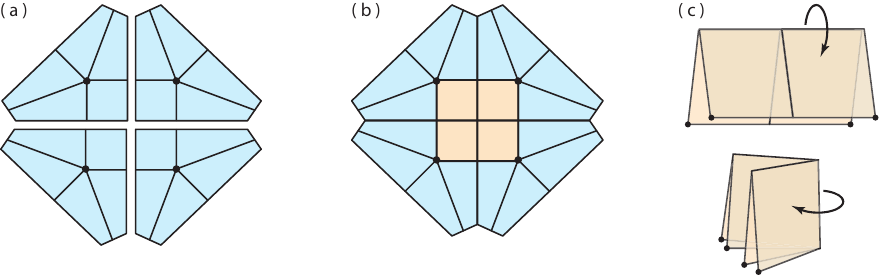}
\caption{Folding cubes of \OM{5}.}
\label{f:map-bhv}
\end{figure}

The small cubes tiling \OM{n+1} are in bijection with \emph{planar} rooted binary trees with labelled leaves, and $\Theta$ can 
be viewed as forgetting planarity conditions.  
The forgetful map $\Theta$ carries out three tasks:  First, it performs an origami folding of the dissected larger square into 
one square, a $2^{n-2} : 1$ folding, as seen in part (c).  Notice that all the centers of each $K_n$ collapse to the unique 
origin of $\iBHV{n}$, the rooted star. 
Second, it identifies this folded cube with its antipodal one, associated to the reflected planar tree, a $2:1$ map.
Finally, the boundaries of the original $K_n$ tiles collapse to the point of $\iBHV{n}$ at $\infty$.
Enumeratively, the folding operation is embodied by the identity
$$n! \ \cdot \ \frac{1}{n}\binom{2n-2}{n-1} \ \cdot \ \frac{1}{2^{n-1}} \ = \ (2n-3)!!$$
taking $n!$ associahedra, partitioning each into $C_{n-1}$ number of cubes, and identify them with a $2^{n-1} : 1$ folding of 
those cubes.
\end{proof}

\begin{figure}[h]
\includegraphics{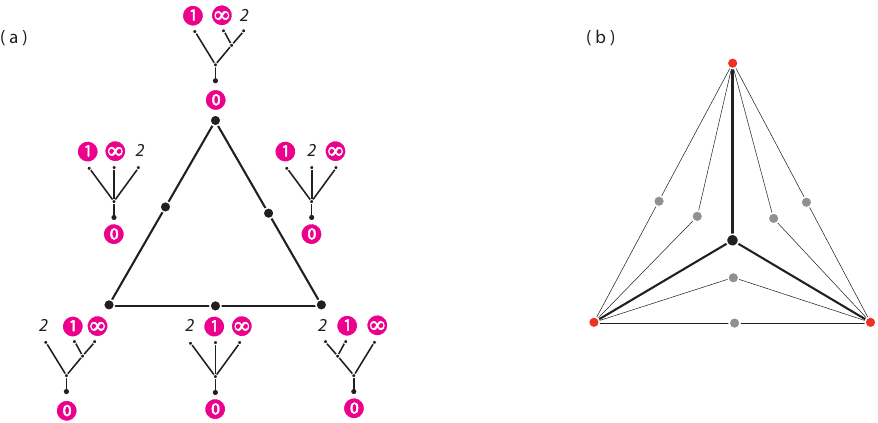}
\caption{The folding of (a) \M{4} into (b) $\iBHV{3}$.}
\label{f:fold}
\end{figure}

\begin{exmp}
Figure~\ref{f:fold}(a) uses the language of rooted trees to redraw \M{4} from Figure~\ref{f:gluing}(b).  Notice how each associahedral $K_3$ edge is decomposed into two cubes (line segments).  Part (b) displays how the six cubes of \M{4} fold into the three edges of $\iBHV{3}$ along each vertex.  Finally, the  centers of the three $K_3$ get identified to the origin of $\iBHV{3}$, and three original vertices of \M{4} glue to form the point at $\infty$.
\end{exmp}

\subsection{}

Two corollaries follow from the Theorem above:

\begin{cor}
There is a natural  map
$$\OM{n+1} \ \to \  \bigvee_{(n-1)!} S^{n-2} $$
folding and collapsing the $n!$ associahedra into $(n-1)!$ spheres.
\end{cor}

\begin{proof}
This follows from Proposition~\ref{p:wedge-sphere} and Theorem~\ref{t:moduli-bhv}.
\end{proof}

\begin{figure}[h]
\includegraphics{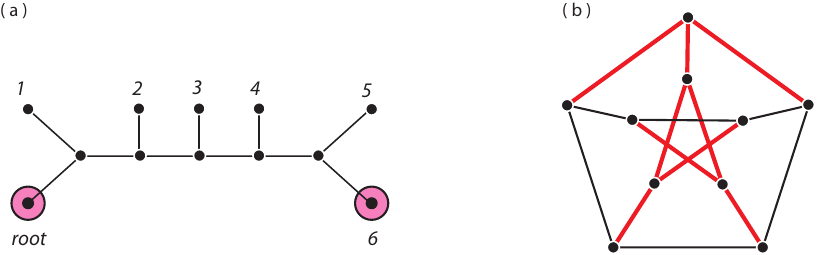}
\caption{(a) Caterpillar trees and (b) contractible chambers of $\T{4}$.}
\label{f:wedge}
\end{figure}

There is an elegant combinatorial description of Eq.~\eqref{e:rw}, based on the work of Robinson and Whitehouse \cite[\S 1]{rw}.  Each chamber of $\iBHV{n}$, being an $(n-2)$-dimensional ball, corresponds to a rooted binary tree with $n$ labeled leaves.   Consider the \emph{caterpillar} tree of Figure~\ref{f:wedge}(a), with one end of the caterpillar\footnote{Detailed relationships between associahedra, permutohedra, and  caterpillars are given in \cite{cdhs}.} as the root, and another end with label $n$ (highlighted leaves of the figure).   There are $(n-1)!$ such chambers of $\iBHV{n}$; the map $\psi$ takes each such chamber to an $(n-2)$-sphere as all other chambers can be contracted to the basepoint.  Figure~\ref{f:wedge}(b) shows the example of $\iBHV{4}$ of Figure~\ref{f:bhv4}(a) where the (red) highlighted edges correspond to noncaterpillars, and thus can be contracted to a point, resulting in Figure~\ref{f:bhv4}(b).

\begin{cor}
The map above factors the $(\Sg_{n-1})$-equivariant composition
\begin{equation} \label{e:comp}
K_n \times_{\Z_{n+1}} \Sg_{n+1} \ \to \ \OM{n+1} \ \to \ \iBHV{n} \ \simeq \ S^1 \wedge
\T{n} \ \sim \ S^{n-2} \wedge \Sg_{n-1}^+ \, .
\end{equation}
\end{cor}

\begin{proof}
This follows from Eq.~\eqref{e:wedge} and Theorem~\ref{t:moduli-bhv}.
\end{proof}

The space $\BHV{n}$ is covered by $(2n-3)!!$ orthants, and the map forgets the order of the pairings in the
Diaconis-Holmes labeling of binary trees. Dividing out the action of $\Sg_{n-1}$ defines a map
\[
K_n \times (\Sg_n/\Sg_{n-1}) \ \to \ S^{n-2} \, ;
\]
this is just an $(\Sg_n/\Sg_{n-1}) = n$-fold copy, parametrized by the cyclic relabelings of the leaves of
trees in $K_n$, of the map
\[
K_n/\partial K_n \ \to \ I^{n-2}/\partial I^{n-2} \ \simeq \ S^{n-2}
\]
which collapses the boundary of the associahedron to a point, and folds the cubical decomposition as
described above.

%
%

\bibliographystyle{amsplain}

\end{document}